\documentclass[11pt]{article}

\usepackage{latexsym}
\usepackage{amsfonts}
\usepackage{amssymb}
\usepackage{mathrsfs}
\usepackage[all]{xy}
\usepackage{epic}
\usepackage{eepic}

\newcommand{\bldw}{\mathbf{w}}

\newcommand{\sgn}{\mbox{sign\,}}

\newcommand{\be}{\begin{equation}}
\newcommand{\ee}{\end{equation}}
\newcommand{\bea}{\begin{eqnarray}}
\newcommand{\eea}{\end{eqnarray}}
\newcommand{\bean}{\begin{eqnarray*}}
\newcommand{\eean}{\end{eqnarray*}}
\newcommand{\brray}{\begin{array}}
\newcommand{\erray}{\end{array}}

\newcommand{\newsection}[1]{\setcounter{equation}{0}
\setcounter{dfn}{0}
\section{#1}}

\newtheorem{dfn}{Definition}[section]
\newtheorem{thm}[dfn]{Theorem}
\newtheorem{lmma}[dfn]{Lemma}
\newtheorem{ppsn}[dfn]{Proposition}
\newtheorem{crlre}[dfn]{Corollary}
\newtheorem{xmpl}[dfn]{Example}
\newtheorem{rmrk}[dfn]{Remark}

\newcommand{\bdfn}{\begin{dfn}\rm}
\newcommand{\bthm}{\begin{thm}}
\newcommand{\blmma}{\begin{lmma}}
\newcommand{\bppsn}{\begin{ppsn}}
\newcommand{\bcrlre}{\begin{crlre}}
\newcommand{\bxmpl}{\begin{xmpl}}
\newcommand{\brmrk}{\begin{rmrk}\rm}

\newcommand{\edfn}{\end{dfn}}
\newcommand{\ethm}{\end{thm}}
\newcommand{\elmma}{\end{lmma}}
\newcommand{\eppsn}{\end{ppsn}}
\newcommand{\ecrlre}{\end{crlre}}
\newcommand{\exmpl}{\end{xmpl}}
\newcommand{\ermrk}{\end{rmrk}}

\newcommand{\bbc}{\mathbb{C}}
\newcommand{\bbz}{\mathbb{Z}}
\newcommand{\bbn}{\mathbb{N}}
\newcommand{\bbr}{\mathbb{R}}

\newcommand{\cla}{\mathcal{A}}

\newcommand{\clh}{\mathcal{H}}

\newcommand{\clk}{\mathcal{K}}
\newcommand{\cll}{\mathcal{L}}

\newcommand{\clg}{\mathcal{G}}

\newcommand{\prf}{\noindent{\it Proof\/}: }

\newcommand{\seq}{\subseteq}

\def \qed { \mbox{}\hfill
$\Box$\vspace{1ex}}

\newcommand{\half}{\frac{1}{2}}

\begin{document}

%%%%%%%%%%%%%%%%%%%%%%%%%%%%%%%%%

%%%%%%%%%%%%%%%%%%%%%%%%%%%%%%%%%
\author{{\sc Partha Sarathi Chakraborty} and
{\sc Arupkumar Pal}}
\title{Torus equivariant spectral triples for odd dimensional quantum spheres
coming from $C^*$-extensions}
\maketitle
%%%%%%%%%%%%%%%%%%%%%%%%%%%%%%%%%%
%%%%%  ABSTRACT
%%%%%%%%%%%%%%%%%%%%%%%%%%%%%%%%%%
 \begin{abstract}
   The torus group $(S^1)^{\ell+1}$ has a canonical action
on the odd dimensional sphere $S_q^{2\ell+1}$.
We take the natural Hilbert space representation
where this action is implemented and characterize
all odd spectral triples acting on that space and equivariant
with respect to that action.
This characterization gives a construction of
an optimum family of equivariant spectral triples having
nontrivial $K$-homology class thus generalizing our earlier results
for $SU_q(2)$. We also relate 
the triple we construct with the $C^*$-extension
\[
 0\longrightarrow \clk\otimes C(S^1)\longrightarrow C(S_q^{2\ell+3})
 \longrightarrow C(S_q^{2\ell+1}) \longrightarrow 0.
\]

 \end{abstract}
{\bf AMS Subject Classification No.:} {\large 58}B{\large 34}, {\large
46}L{\large 87}, {\large
  19}K{\large 33}

%%%%%%%%%%%%%%%%%%%%%%%%%%%%%%%%%%%%%%%%%%%%%%%%%%%

% \tableofcontents

%%%%%%%%%%%%%%%%%%%%%%%%%%%%%%%%%%%%%%%%%%%%%%%%%%%
\newsection{Introduction}
In noncommutative geometry (NCG),  a geometric space is
described by a triple $(\cla,\clh,D)$, called a spectral triple, 
with $\cla$ being an involutive
algebra represented as bounded operators on a Hilbert space $\clh$,
and $D$ being a selfadjoint operator with compact resolvent and having
bounded commutators with the algebra elements. The operator $D$ should be
nontrivial in the sense that the associated Kasparov module should
give a nontrivial element in $K$-homology. A natural question is, 
are there enough spectral triples around us? The answer is 
both yes and no. If we do not demand any further properties 
then by a theorem of Baaj and Julg (\cite{b-j}), given any 
countable subalgebra $\cla$ of a $C^*$-algebra there exists a 
spectral triple $(\cla, \clh, D)$. But if we demand further 
properties like finite summability then given a dense 
subalgebra of a $C^*$-algebra it may not admit a finitely 
summable spectral triple (\cite{co0}). Therefore 
given a natural dense subalgebra of a $C^*$-algebra it is 
meaningful to ask whether it admits finitely summable nontrivial 
spectral triples. 
Also,  the result of Baaj \& Julg starts from a Fredholm module,
so one has very little control over the Hilbert space
or the representation.

In an earlier paper (\cite{c-p3}), the authors studied spectral triples
for the odd dimensional quantum spheres taking the Hilbert space
to be the $L_2$ space of the sphere and the representation
to be the natural representation by left multiplication there.
In the present article, we fix a different representation space
dictated by the torus action on the sphere, and investigate
spectral triples for that. The results here generalize those in~\cite{c-p2}.

We will use the method described in \cite{c-p3} and used implicitly
in \cite{c-p1} and \cite{c-p2}.
Observe that the self-adjoint operator $D$ in a spectral
triple comes with two very crucial restrictions on it, namely, it has to have
compact resolvent, and must have bounded commutators with algebra
elements. Various analytic consequences of the compact resolvent
condition (growth properties of the commutators of the algebra
elements with the sign of $D$) have been used in the past by various
authors.   We will  exploit
it from a combinatorial point of view.  The idea is very simple. Given a
selfadjoint operator with compact resolvent, one can associate with it
a certain  graph in a natural way.  This makes it possible to do a detailed
combinatorial analysis of the growth restrictions (on the eigenvalues
of $D$) that come from the boundedness of the commutators, and to
characterize the sign of the operator $D$ completely.

 We take a representation space where the canonical action
of $(S^1)^{\ell+1}$ on $C(S_q^{2\ell+1})$ is implemented. 
If we further want our Dirac operator $D$ to be equivariant 
with respect to the torus action then $D$ should commute 
with the unitaries implementing that action. Hence $D$ 
respects the spectral subspaces. This allows us to write 
down the form of the Dirac operator. Then using the 
boundedness of the commutators we 
 completely characterize all equivariant Dirac operators.
We also produce a nontrivial optimal equivariant Dirac.

Odd dimensional quantum spheres of successive dimension
are related through a short exact sequence that says that
the $(2\ell+3)$-dimensional sphere $C(S_q^{2\ell+3})$ is an
extension of the $(2\ell+1)$-dimensional sphere $C(S_q^{2\ell+1})$
by $C(S^1)$. One can naturally associate a $KK_1(C(S_q^{2\ell+1}),C(S^1))$ element with
such an extension. In the last section, we compute this KK-element
and show that the generic spectral triple that we construct in section~3
comes from this KK-element.

%%%%%%%%%%%%%%%%%%%%%%%%%%%

%%%%%%%%%%%%%%%%%%%%%%%%%%%%%%%%%%%%%%%%%%%%%%%%%%%%%%
\newsection{Torus action on quantum spheres}
%%%%%%%%%%%%%%%%%%%%%%%%%%%%%%%%%%%%%%%%%%%%%%%%%%%%%%
Let $q\in[0,1]$.
The $C^*$-algebra $A_\ell=C(S_q^{2\ell+1})$ of continuous functions
on the quantum sphere $S_q^{2\ell+1}$
is the universal $C^*$-algebra generated by
elements
$z_1, z_2,\ldots, z_{\ell+1}$
satisfying the following relations (see~\cite{h-s}, \cite{v-s}):
\bea \label{Reln}
z_i z_j & =& qz_j z_i,\qquad 1\leq j<i\leq \ell+1, \nonumber \\
z_i z^{*}_j & =& q z^{*}_j z_i ,\qquad 1\leq i\neq j\leq \ell+1, \nonumber\\
z_i z^{*}_i - z^{*}_i z_i +
(1-q^{2})\sum_{k>i} z_k z^{*}_k &=& 0,\qquad \hspace{2em}1\leq i\leq \ell+1,\\
\sum_{i=1}^{\ell+1} z_i z^{*}_i &=& 1. \nonumber
\eea
Let $N$ be the number operator given by $N: e_n \mapsto n e_n$ on $L_2(\bbn)$ 
and $S$ be the  shift $S: e_n \mapsto e_{n-1}$. 
We will use the same symbol $S$ to denote shift 
on $L_2(\bbn)$ as well as on $L_2(\bbz)$. In the case 
of $L_2(\bbn)$, $S(e_0)$ is defined to be zero. 
Let 
\[
 \clh_\ell=\underbrace{L_2(\bbn)\otimes\cdots\otimes
         L_2(\bbn)}_{\ell \mbox{ copies}}\otimes L_2(\bbz).
\]
Let $\pi_\ell$ be the representation of $A_\ell$ on the space $\cll(\clh_\ell)$
of bounded operators on $\clh_\ell$ given on the generators by
\bean
 z_k &\mapsto & \underbrace{q^N\otimes\ldots\otimes q^N}_{k-1 \mbox{ copies}}\otimes
      \sqrt{1-q^{2N}}S^*\otimes 
   \underbrace{I \otimes\cdots\otimes I}_{\ell+1-k \mbox{ copies}},\quad
            1\leq k\leq \ell,\\
z_{\ell+1} &\mapsto &
    \underbrace{q^N\otimes\cdots\otimes q^N}_{\ell \mbox{ copies}}
       \otimes S^*,
\eean
Then $\pi_\ell$ gives a faithful representation
of $A_\ell$ on $\clh_\ell$ (see lemma~4.1 and remark~4.5, \cite{h-s}).
Observe that for all $a\in A_\ell$, the operators $\pi_\ell(a)$ actually
lift to operators on the Hilbert $C(S^1)$-module $L_2(\bbn)\otimes C(S^1)$.

$K$-groups of these $C^*$-algebras have been computed by
Vaksman \& Soibelman and Hong \& Szymanski:
\bppsn[\cite{v-s},\cite{h-s}]\label{ktheory}
$K_0(A_\ell)=K_1(A_\ell)=\bbz$.
\eppsn

The group $(S^1)^{\ell+1}$ has an action on $C(S_q^{2\ell+1})$ given on the generating
elements by
\[
 \tau_\bldw(z_i)=w_i z_i,\qquad
\bldw=(w_1,w_2,\ldots,w_{\ell+1})\in(S^1)^{\ell+1}.
\]
If $U_\bldw$ denotes the unitary
$w_1^N\otimes w_2^N\otimes\cdots \otimes w_{\ell+1}^N$
on $\clh_\ell$, then one has $\pi_\ell(\tau_\bldw(a))=U_\bldw \pi_\ell(a) U_\bldw^*$
for all $a\in C(S_q^{2\ell+1})$.
Thus $(\pi_\ell, U)$ is a covariant representation of 
$(A_\ell,(S^1)^{\ell+1},\tau)$ on $\clh_\ell$.
In the next section, we characterize all equivariant spectral
triples for this representation and construct an optimal
triple using this characterization.

%%%%%%%%%%%%%%%%%%%%%%%%%%%%%%%%%%%%%%%%%%
\section{Equivariant spectral triples}
%%%%%%%%%%%%%%%%%%%%%%%%%%%%%%%%%%%%%%%%%%
Let $\Gamma=\underbrace{\bbn\times\cdots\times\bbn}_{\ell\mbox{ copies}}\times\bbz$,
so that $L_2(\Gamma)=\clh_\ell$.
For $\gamma=(\gamma(1),\gamma(2),\cdots,\gamma(\ell+1))\in\Gamma$,  $e_\gamma$ denotes the basis element of $\clh_\ell$ given by
$e_{\gamma(1)}\otimes\cdots\otimes e_{\gamma(\ell+1)}$.

\bthm
Let $D$ be a self-adjoint operator with compact resolvent
on $\clh_\ell$ that commutes with the operators $U_\bldw$.
Then $D$ must diagonalise with respect to
the canonical basis, i.\ e.\ must be of the form
\be
e_\gamma\mapsto
  d(\gamma)e_\gamma,
\ee
where $d(\gamma)\in\bbr$ for all $\gamma\in\Gamma$.

Moreover, such an operator $D$ will have bounded commutators
with elements from the *-subalgebra of $C(S_q^{2\ell+1})$
generated by the $z_i$'s if and only if the $d(\gamma)$'s
obey the following condition:
\be\label{eqbdd_sph_neq}
| d(\gamma)- d(\gamma+\epsilon_k)| = O(q^{-\gamma(1)-\ldots - \gamma(k-1)}),
      \qquad 1\leq k\leq \ell+1,
\ee
where $\epsilon_k$ stands for the vector whose k\raisebox{.4ex}{th}
coordinate is 1 and all other coordinates are 0.
\ethm
\prf
The first part is immediate. For the second part,
just observe that
\bean
[D,\pi(z_k)]e_\gamma &=&(d(\gamma+\epsilon_k)-d(\gamma))
       q^{\gamma(1)+\ldots + \gamma(k-1)}\sqrt{1-q^{2\gamma(k)+2}}
        e_{\gamma+\epsilon_k},\quad 1\leq k\leq \ell, \cr
[D,\pi(z_{\ell+1})]e_\gamma &=&   (d(\gamma+\epsilon_{\ell+1})-d(\gamma))
       q^{\gamma(1)+\ldots + \gamma(\ell)}  e_{\gamma+\epsilon_{\ell+1}}.
\eean
\qed\\
By a compact perturbation, one can ensure that all
the $d(\gamma)$'s are  nonzero in the above theorem.
We will assume from now on that $d(\gamma)\neq 0$ for all $\gamma$.
Using~(\ref{eqbdd_sph_neq}) we get a constant $c$ such that 
$| d(\gamma)- d(\gamma+\epsilon_k)|q^{-\gamma(1)-\ldots - \gamma(k-1)} < c$,
with $\epsilon_k$ as in the theorem. Now join two elements
$\gamma$ and $\gamma'$ in $\Gamma$ by an edge if
$|d(\gamma)-d(\gamma')|\leq c$. Call the resulting graph $\clg$ the
growth graph for $D$.

\blmma\label{p_path1}
%%%%%%%%%%%%
Let $k$ be an integer with  $1 \leq k\leq  \ell+1$.
Let
\[
\gamma=(0,\ldots,0,r,i_{k+1},\ldots,i_{\ell+1}),\quad
\gamma'=(0,\ldots,0,s,i_{k+1},\ldots,i_{\ell+1}).
\]
Then there is a path in $\clg$ of length $|r-s|$ joining $\gamma$ and
$\gamma'$ such that all vertices on this path are of the form
$(0,\ldots,0,t,i_{k+1},\ldots,i_{\ell+1})$.
%%%%%%%%%%%%
\elmma
\prf
Assume without loss in generality that
$\gamma(k)<\gamma'(k)$. Write $r=\gamma'(k)-\gamma(k)$.
From~(\ref{eqbdd_sph_neq}), it is clear that
if $\delta(i)=0$ for $1\leq i\leq k-1$, then
there is an edge joining $\delta$ and $\delta+\epsilon_k$.
Thus
$(\gamma,\gamma+\epsilon_k, \gamma+2\epsilon_k,\ldots,\gamma+r\epsilon_k)$
will give us a required path.
\qed

\blmma \label{p_path2}
%%%%%%%%%%%%%%%
Let $k$ be an integer with  $1 \leq k\leq  \ell+2$.
Let
\[
\gamma=(i_1,\ldots,i_{k-1},i_{k},\ldots,i_{\ell+1}),\quad
\gamma'=(0,\ldots,0,i_{k},\ldots,i_{\ell+1}).
\]
Then there is a path of length $|i_1|+\ldots+|i_{k-1}|$ joining $\gamma$ and
$\gamma'$ such that all vertices on this path are of the form
$(j_1,\ldots,j_{k-1},i_{k},\ldots,i_{\ell+1})$, where
each $j_n$ lies between 0 and $|i_n|$.
%%%%%%%%%%%%%%%
\elmma
\prf
For $1\leq j\leq k$, let $\gamma_j$ denote
the element of $\Gamma$ whose
first $j-1$ coordinates are 0 and $j$\raisebox{.4ex}{th}
coordinate onwards coincide with those of $\gamma$.
Thus $\gamma_1=\gamma$ and $\gamma_k=\gamma'$.
Now apply the previous proposition to get
a path of length $|\gamma_j(j)-\gamma_{j+1}(j)|=\gamma(j)$
joining $\gamma_j$ and $\gamma_{j+1}$ for $1\leq j\leq k-1$.
Joining all these paths together, one gets the required path.
\qed

\bppsn\label{optimal1}
Let $D$ be a Dirac operator that commutes with the operators
$U_\bldw$. Then $D$ must be of the form
$e_\gamma\mapsto d(\gamma)e_\gamma$
where
\[
|d(\gamma)|=O(\gamma(1)+\ldots +\gamma(\ell)+|\gamma(\ell+1)|+1)|.
\]
\eppsn
\prf
Note that if $\gamma$ is an arbitrary element of the growth graph $\clg$, then by the previous lemmas $\gamma$ can be connected with $0$ by a path of length $ \gamma(1)+\ldots +\gamma(\ell)+|\gamma(\ell+1)$, hence the result.
\qed

\bthm\label{sign_sph}
Write $\Gamma^+=\{\gamma\in\Gamma: d(\gamma)>0\}$,
and $\Gamma^-=\Gamma\verb1\1\Gamma^+$.
There exist nonnegative integers $M_1,M_2,\ldots,M_{\ell+1}$
such that for each $k\in\{1,2,\ldots,\ell\}$ and for each
\[
(i_{k+1},i_{k+2},\ldots,i_{\ell+1})\in F_k:=
    \prod_{r=k+1}^\ell\{0,1,\ldots,M_r\}\times
     \{-M_{\ell+1},-M_{\ell+1}+1,\ldots,M_{\ell+1}\},
\]
none of the following sets intersect both
$\Gamma^+$ and $\Gamma^-$:
\[
A_1=\{\gamma\in\Gamma: \gamma(\ell+1)> M_{\ell+1}\},
\quad
A_2=\{\gamma\in\Gamma: \gamma(\ell+1)< -M_{\ell+1}\},
\]
\[
B_{k,(i_{k+1},i_{k+2},\ldots,i_{\ell+1})}=
      \{\gamma\in\Gamma: \gamma(k)> M_k, \gamma(r)=i_r
   \mbox{ for }k+1\leq r\leq \ell+1\}.
\]
\ethm
\prf
We will construct these numbers $M_1,M_2,\cdots M_{\ell+1}$ inductively starting from $M_{\ell+1}$.
Assume if possible there are two sequences of elements
$\gamma_k\in\Gamma^+$ and $\delta_k\in\Gamma^-$ such that
\[
% \cdots \gamma_{-1}(1)< \delta_{-1}(1) <
 \gamma_0(\ell+1)
     < \delta_0(\ell+1) < \gamma_1(\ell+1) < \delta_1(\ell+1) < \cdots.
\]
For each $k$, use lemma~\ref{p_path2} to get a path $p_k$
from $\gamma_k$ to $\delta_k$ such that for any vertex on the path,
the $\ell+1$\raisebox{.4ex}{th} coordinate lies between $\gamma_k(\ell+1)$ and $\delta_k(\ell+1)$.
This last condition would ensure that the paths $p_k$ are all disjoint.
Since $p_k$ connects points of $\Gamma^+$ with $\Gamma^-$, there is a vertex $\mu_k$ in $p_k$ such that $d(\mu_k) \in [-c,c]$. Moreover disjointness of the $p_k$'s implies that the vertices $\mu_k$ are all distinct. Therefore counted with multiplicity, the compact interval $[-c,c]$ has infinitely many eigenvalues of $D$, a contradiction to compact resolvent condition for $D$. Therefore there exists $M^\prime_{\ell+1} $ such that $\{\gamma\in\Gamma: \gamma(\ell+1)> M^\prime_{\ell+1}\}$ does not intersect both $\Gamma^+$ and $\Gamma^-$.  
One can similarly show that if there are elements
$\gamma_k\in\Gamma^+$ and $\delta_k\in\Gamma^-$ such that
\[
 \gamma_0(\ell+1)
     > \delta_0(\ell+1) > \gamma_1(\ell+1) > \delta_1(\ell+1) > \cdots,
\]
 then  there is some big enough natural number $M^{\prime\prime}_{\ell+1}$
such that the set $\{\gamma\in\Gamma: \gamma(\ell+1)< -M^{\prime \prime}_{\ell+1}\}  $ is either in $\Gamma^+$ or in $\Gamma^-$. Now taking $M_{\ell+1}=\max\{M^\prime_{\ell+1}, M^{\prime \prime}_{\ell+1} \}$, we get that neither of $A_1,A_2$ intersect both $\Gamma^+$ and $\Gamma^-$.

Next, given $M_{k+1},\ldots, M_{\ell+1}$
and
$(i_{k+1},i_{k+2},\ldots,i_{\ell+1})\in F_k$,
if there are elements $\gamma_n\in\Gamma^+$ and $\delta_n\in\Gamma^-$
with
\[
\gamma_n(j)=i_j=\delta_n(j), \quad k+1\leq j\leq \ell+1,
\]
\[
\gamma_0(k)<\delta_0(k) <\gamma_1(k) <\delta_1(k)<\cdots,
\]
then using lemma~\ref{p_path2} again, one can join each pair
 $(\gamma_n,\delta_n)$ by disjoint paths and  arguing as above  arrive at a contradiction to the fact that $D$ has compact resolvent. 
Therefore the existence of $M_k$ follows.
\qed

\bthm
Let $D_{torus}$ be the operator $e_\gamma\mapsto d(\gamma)e_\gamma$ on $\clh_\ell$
where the $d(\gamma)$'s are given by
\[
d(\gamma)=\cases{\gamma(1)+\ldots +\gamma(\ell)+|\gamma(\ell+1)|
                        & if $\gamma(\ell+1)\geq 0$,\cr
          -( \gamma(1)+\ldots +\gamma(\ell)+|\gamma(\ell+1)|)
                & if      $\gamma(\ell+1) < 0$.}
\]
Then $(C(S_q^{2\ell+1}),\clh_\ell, D_{torus})$ is a
nontrivial $(\ell+1)$-summable spectral triple.

The operator $D_{torus}$ is optimal, i.\ e.\ if $D$ is
any Dirac operator acting on $\clh$ that
commutes with the $U_\bldw$'s, then there exist positive
reals $a$ and $b$ such that
\[
|D|\leq a+b|D_{torus}|.
\]
\ethm
\prf
Clearly $D_{torus}$ is a selfadjoint operator with compact resolvent.
That it has bounded commutators with the $\pi(z_j)$'s
follow by direct verification.

From the commutation relations that the generators $z_j$ obey,
it follows that $z_{\ell+1}$ is normal and the element
 $z_{\ell+1}^* z_{\ell+1}$ has spectrum $\{q^{2n}:n\in\bbn\}\cup\{0\}$.
Let
\[
u=\chi_{\{1\}}(z_{\ell+1}^* z_{\ell+1})(z_{\ell+1}-1) +1.
\]
It is easy to see that $u$ is a unitary. We will now compute the pairing
between $D_{torus}$ and $\pi(u)$.
First observe that the action of $\pi(u)$ on $\clh$ is given by
\[
\pi(u)e_\gamma=\cases{e_{\gamma+\epsilon_{\ell+1}}
          & if $\gamma(i)= 0$ for  $1\leq i\leq\ell$,\cr
          e_\gamma & otherwise.}
\]
Write $P=\half(I+\sgn D_{torus})$. Then $P$ is the projection onto
the closed linear span of $\{e_\gamma: \gamma(\ell+1)\geq 0\}$.
It follows that the index of $PuP$ is $-1$.

Summability follows from the observation that
the number of elements in
$\{(i_1,\ldots,i_{\ell+1})\in \bbn^\ell\times\bbz:
    \sum_{k=1}^\ell i_k +|i_{\ell+1}|\leq n\}$ is of the order $n^{\ell+1}$.

Optimality is a consequence of proposition~\ref{optimal1}.
\qed

\bthm
Let $D$ be a Dirac operator on $\clh$ that commutes with the operators
$U_\bldw$. Then either $D$ is trivial or
has the same $K$-homology class as $D_{torus}$ or $-D_{torus}$.
\ethm
\prf
If $D$ is a self-adjoint operator with compact resolvent
on $\clh$ that commutes with the operators $U_\bldw$ and if
$P=\half(\sgn D+I)$, then by theorem~\ref{sign_sph},
$P$ is the projection onto the closed linear span
of $\{e_\gamma: \gamma\in \Gamma^+\}$ where $\Gamma^+$
must be of one of the following form:
\be
A_1\cup \bigl(\cup_{x\in E}B_x\bigr),
\ee
\be
A_2\cup (\cup_{x\in E}B_x),
\ee
\be
A_1\cup A_2\cup \left(\cup_{x\in E}B_x\right),
\ee
\be
\cup_{x\in E}B_x,
\ee
where $E$ is some finite subset of $\cup_{k=1}^\ell \{k\}\times F_k$.
By direct calculations in the first two cases  the index of $P\pi(u)P$ turns out to be $-1$ and $1$ respectively,
whereas in the last two cases, the index is zero.
Thus one always has
\[
\langle [u],(C(S_q^{2\ell+1}),\clh, D)\rangle = 0 \mbox{ or }\pm 1.
\]
By~\cite{r-s}, we have $K^1(C(S_q^{2\ell+1}))=\bbz$. therefore
the result follows.
\qed

%%%%%%%%%%%%%%%%%%%%%%%%%%%%%%%%%%%%%%%%%%%%%%%%%%%%%
%%%%%%%%%%%%%%%%%%%%%%%%%%%%%%%%%%%%%%%%%%%%%%%%%%%%%
\section{Relation with $C^*$-extensions}
%%%%%%%%%%%%%%%%%%%%%%%%%%%%%%%%%%%%%%%%%%%%%%%%%%%%%
In this section we will denote the generators  for $A_\ell$
by $z_k$ and the generators for $A_{\ell+1}$ by $y_k$. 
$A_\ell^0$ will denote the *-subalgebra of $A_\ell$
generated by the $z_k$'s.
Let $J_\ell^0$ denote the two-sided *-ideal in $A_\ell^0$
generated by $z_{\ell+1}$ and let $J_\ell$ denote
the norm closure of $J_\ell^0$ in $A_\ell$. Thus $J_\ell$
is the ideal in $A_\ell$ generated by the element $z_{\ell+1}$.

For a Hilbert $C^*$-module $E$, we will denote by $\cll(E)$
the $C^*$-algebra of bounded adjointable operators on $E$, and by
$\clk(E)$ its ideal of `compact' operators. We denote by $\clk$
the $C^*$-algebra $\clk(\clh)$ for an infinite dimensional Hilbert space $\clh$.

\blmma \label{ker}
Let $C^*(S)$ denote the $C^*$-algebra generated by the
operator $S$ on $L_2(\bbz)$.
Then one has $J_\ell\cong \clk(L_2(\bbn^\ell)) \otimes C^*(S)\cong \clk\otimes C(S^1)$.
\elmma
\prf
We will identify $A_\ell$ with $\pi_\ell(A_\ell)$.

For $1\leq k\leq \ell$, denote by $X_k$ the operator
\[
\underbrace{q^N\otimes\ldots\otimes q^N}_{\mbox{$k$ copies}}\otimes 
\underbrace{I\otimes\ldots\otimes I}_{\mbox{$\ell+1-k$ copies}}
\]
on $\clh_\ell$.
Write $X_0=I$. Then it is easy to check that one
has the relations
\[
 z_k z_k^*=X_{k-1}^2-X_k^2,\quad 1\leq k\leq \ell.
\]
It follows that  $X_k\in A_\ell$ for all $1\leq k\leq \ell$.

Write $p_{ij}$ for the rank one operator $|e_i\rangle\langle e_j|$
on $L_2(\bbn)$.
Then
\[
 p_{i_1j_1}\otimes\ldots \otimes p_{i_\ell j_\ell}\otimes S^k
\]
can be written in the form
\[
f_1(X_1)\ldots f_\ell(X_\ell)z_{\ell+1}^{-k}g_1(X_1)\ldots g_\ell(X_\ell)
\]
where $f_i$, $g_i$ are continuous functions on the spectrums of the respective
$X_i$'s.
Therefore 
$p_{i_1j_1}\otimes\ldots \otimes p_{i_\ell j_\ell}\otimes S^k \in J_\ell$.
It follows from this that
$\clk(L_2(\bbn^\ell))\otimes C^*(S)\seq J_\ell$.

For the reverse inclusion, observe that any polynomial
in the $z_i$'s and their adjoints is a finite sum of the form
$\sum_j T_j\otimes S^{k_j}$ where $T_j\in\cll(L_2(\bbn^\ell))$ and $k_j\in \bbz$.
Therefore $J_\ell^0$ is contained in $\clk(L_2(\bbn^\ell))\otimes C^*(S)$.
Same is therefore true for its closure $J_\ell$.
\qed

\bppsn \label{shexseq}
Let $\sigma_\ell:\cla_{\ell+1} \rightarrow \cla_\ell$ be the homomorphism given by
 \[
 y_i\mapsto \cases{z_i & if $1\leq i\leq \ell+1$,\cr
                0 & if $i=\ell+2$.}
\]
Then we have the following short exact sequence
\be\label{eq:shexseq}
 0\longrightarrow J_{\ell+1}\longrightarrow A_{\ell+1}
  \stackrel{\sigma_\ell}{\longrightarrow} A_\ell\longrightarrow 0.
\ee
\eppsn

We will need the following lemma for the proof.

\blmma \label{quotient}
Let $\cla$ be the universal $C^*$-algebra in 
noncommuting variables $x_1,x_2,\cdots x_n$ subject 
to relations $R_1(x_1,x_2,\cdots,x_n),\cdots, R_j(x_1,x_2,\cdots,x_n)$. 
Let $J$ be the ideal of $\cla$ generated by  noncommutative polynomials $Q_1(x_1,x_2,\cdots,x_n),Q_2(x_1,x_2,\cdots,x_n),\cdots, Q_k(x_1,x_2,\cdots,x_n)$. 
Then $\cla/J$ is isomorphic to the universal $C^*$-algebra $\cla(J)$ 
generated by $x_1,x_2,\cdots,x_n$ subject to the relations 
$R_1,\cdots,R_j,Q_1,\cdots,Q_k$. 
\elmma
Note that it is part of the 
hypothesis that the universal $C^*$-algebras $\cla$ and $\cla(J)$ exist.\\[1ex]
\prf
Let  $\xi_1,\cdots,\xi_n$ be the generating elements of 
$\cla(J)$. Clearly we have a surjection $q : \cla(J) \rightarrow \cla/J $ 
mapping $\xi_i$ to $x_i$. To show that this is injective 
it is enough to show that given a polynomial 
$\alpha=f(\xi_1,\cdots,\xi_n) \in \cla(J)$, 
one has  $\|q(\alpha)\|=\|a\|,$ 
where $a=f(x_1,\cdots,x_n)$. 
Now observe that
\bean
\|a\|& =&\sup\{\|\pi(a)\|: \pi \mbox{ is a representation of } \cla, \pi(J)=0 \}\\
& = & \sup \{ \|\pi(a)\|: \pi \mbox{ is a representation of the algebra generated by } x_1,x_2,\cdots x_n \\ 
&  &\mbox{ subject to  }
 R_1,\cdots,R_j,Q_1,\dots , Q_k \}\\
 & = & \|\alpha\|.
\eean
Thus the proof is complete.
\qed

\noindent \textit{Proof of proposition~\ref{shexseq}.}
Clearly $J_{\ell+1} \subseteq \ker(\sigma_\ell)$ and lemma~\ref{quotient} gives $\cla_{\ell+1}/J_{\ell+1} \cong \cla_{\ell+1}(J_{\ell+1})$. Also note that
in the defining relations for the generators for $\cla_{\ell+1}$ if we put $y_{\ell+2}=0$ we get the relations for $\cla_\ell$, hence $ \cla_{\ell+1}(J_{\ell+1})=\cla_\ell$.  
Therefore $\ker(\sigma_\ell)=J_{\ell+1}$, hence the result.
\qed

Proposition~\ref{shexseq} gives a homomorphism 
$\psi_{\ell+1}:\cla_{\ell+1}\rightarrow M(J_{\ell+1})$. 
Using lemma~\ref{ker} we get $M(J_{\ell+1})\cong \cll(L_2(\bbn^{\ell+1}) \otimes C(S^1))$. 
Thus $\psi_{\ell+1}$ is  given by:
\bean
y_k &\mapsto& \underbrace{q^N\otimes\cdots\otimes q^N}_{k-1 \mbox{ copies}}\otimes
 \sqrt{1-q^{2N}}S^*\otimes  
\underbrace{I \otimes\cdots\otimes I}_{\ell+2-k \mbox{ copies}},\quad
            1\leq k\leq \ell+1,\cr
y_{\ell+2} &\mapsto&
    \underbrace{q^N\otimes\cdots\otimes q^N}_{\ell+1 \mbox{ copies}}
       \otimes Z.
\eean
Here $Z:C(S^1) \rightarrow C(S^1)$ denotes the operator given by $(Zf)(z)=zf(z)$.

Define $\tilde{\sigma}_\ell:A_\ell\rightarrow \cll(\clh_\ell\otimes C(S^1))$
by
\bean
z_k & \mapsto & \underbrace{q^N\otimes\cdots\otimes q^N}_{k-1 \mbox{ copies}}\otimes
      \sqrt{1-q^{2N}}S^*\otimes 
     \underbrace{I \otimes\cdots\otimes I}_{\ell+2-k \mbox{ copies}},\quad 1\leq k\leq \ell,\\
z_{\ell+1} & \mapsto & \underbrace{q^N\otimes\cdots\otimes q^N}_{\ell \mbox{ copies}}\otimes
      S^*\otimes I.
\eean
Let
\[
 E_\ell =\underbrace{\clk(L_2(\bbn))\otimes\cdots\otimes \clk(L_2(\bbn))}_{\ell \mbox{ copies}}
   \otimes C(S^1), \quad
F_\ell =\underbrace{L_2(\bbn)\otimes\cdots\otimes L_2(\bbn)}_{\ell \mbox{ copies}}
     \otimes C(S^1).
\]
Let $U$ be the unitary from $L_2(\bbn)\oplus L_2(\bbn)$ onto $L_2(\bbz)$
given by
\[
e_n\oplus 0 \mapsto e_n,\qquad 0\oplus e_n\mapsto e_{-n-1}, \qquad n\in\bbn.
\]
Using this unitary in the $(\ell+1)$\raisebox{.4ex}{th} copy,
one can identify $\clh_\ell\otimes C(S^1)$ with $F_{\ell+1}\oplus F_{\ell+1}$
Let $P\in\cll(L_2(\bbz))$ be the projection onto the $L_2(\bbn)$ part
and let $Q_\ell=\underbrace{I\otimes\cdots\otimes I}_{\ell \mbox{ copies}}\otimes P\otimes I$. 
Define $C_\ell:\cll(\clh_\ell\otimes C(S^1))\rightarrow \cll(F_{\ell+1})$
by $C_\ell(T)=Q_\ell T Q_\ell$.
Now define $\hat{\sigma}_\ell: A_\ell\rightarrow \cll(F_{\ell+1})$
by $\hat{\sigma}_\ell(a)=C_\ell \tilde{\sigma}_\ell(a)$.
For convenience, we summarize various maps and the spaces between which
they act in the following diagram:

\hspace*{4em}
\def\labelstyle{\scriptstyle}
\xymatrix@C=25pt@R=40pt{
 & J_{\ell+1}\ar[r]\ar@{<->}[d]^{=} & A_{\ell+1}\ar[r]^{\sigma_\ell}\ar@{<->}[d]^{\psi_{\ell+1}} & 
         A_\ell \ar[ddl]_{\hat{\sigma}_\ell} \ar[d]^{\tilde{\sigma}_\ell}&  \\
&E_{\ell+1}\ar[d]^{\seq}&  \psi_{\ell+1}(A_{\ell+1})\ar[d]^{\seq} & \cll(\clh_\ell\otimes C(S^1))\ar@{<->}[d]^{\cong}& \\
& M(E_{\ell+1}) \ar@{<->}[r]^{\cong} & \cll(F_{\ell+1}) &  \cll(F_{\ell+1}\oplus F_{\ell+1})\ar[l]_{C_\ell}
}

\bthm
The element $(\clh_\ell\otimes C(S^1), \tilde{\sigma}, 2Q-I)$
gives the $KK$-class in $KK^1(C(S_q^{2\ell+1}),C(S^1))$
corresponding to the extension~(\ref{eq:shexseq}).
\ethm
\prf
Let $r\in\bbn$ and let $p$ be a polynomial in
noncommuting variables and their adjoints.
Using the observation that $Q_\ell$ commutes with
$\tilde{\sigma}_\ell(z_k)$ for $1\leq k\leq \ell$,
one gets
\begin{enumerate}
\item $\hat{\sigma}_\ell(z_{\ell+1}^r p(z_1,\cdots,z_\ell,z_1^*,\cdots,z_\ell^*))=
    \hat{\sigma}_\ell(z_{\ell+1}^r)\hat{\sigma}_\ell(p(z_1,\cdots,z_\ell,z_1^*,\cdots,z_\ell^*)).$
\item $\hat{\sigma}_\ell({(z_{\ell+1}^*)}^r p(z_1,\cdots,z_\ell,z_1^*,\cdots,z_\ell^*))=
  \hat{\sigma}_\ell({(z_{\ell+1}^*)}^r)\hat{\sigma}_\ell(p(z_1,\cdots,z_\ell,z_1^*,\cdots,z_\ell^*)).$
\end{enumerate}
Using this one can now easily show that
\begin{enumerate}
\item
$ \hat{\sigma}_\ell(p(z_1,\cdots,z_\ell,z_1^*,\cdots,z_\ell^*)) =
   \psi_{\ell+1}(p(y_1,\cdots,y_\ell,y_1^*,\cdots,y_\ell^*)).$
\item $
\hat{\sigma}_\ell(z_{\ell+1}^r)-\psi_{\ell+1}(y_{\ell+1}^r) \in \clk(L_2(\bbn^{\ell+1}))\otimes C^*(S)=\psi_{\ell+1}(J_{\ell+1}).$
\item 
$ \hat{\sigma}_\ell({(z_{\ell+1}^*)}^r)-\psi_{\ell+1}({(y_{\ell+1}^*)}^r) \in \clk(L_2(\bbn^{\ell+1}))\otimes C^*(S)=\psi_{\ell+1}(J_{\ell+1}).$
\end{enumerate}
It follows from these that for any polynomial $p$  we have
\bea \label{lift}
\lefteqn{\hspace*{-10em}\hat{\sigma}_\ell(p(z_1,\cdots,z_{\ell+1},z_1^*,\cdots,z_{\ell+1}^*))
 -\psi_{\ell+1}(p(y_1,\cdots,y_{\ell+1},y_1^*,\cdots,y_{\ell+1}^*)) }\nonumber\\
&\in& \clk(L_2(\bbn^{\ell+1}))\otimes C(S^1)=\psi_{\ell+1}(J_{\ell+1}).
\eea

Let $\tau:\cla_{\ell} \rightarrow M(J_{\ell+1})/J_{\ell+1}$ be the 
Busby invariant for the extension (\ref{eq:shexseq}),
and let $\Phi: M(J_{\ell+1})\rightarrow M(J_{\ell+1})/J_{\ell+1}$ be the quotient map. 
For a polynomial $p$ in noncommuting variables and their adjoints, 
we now have from~(\ref{lift}),
 \bean
 \tau(p(z_1,\cdots,z_{\ell+1},z_1^*,\cdots,z_{\ell+1}^*)) 
    & = & \Phi \circ \psi_\ell(p(y_1,\cdots,y_{\ell+1},y_1^*,\cdots,y_{\ell+1}^*)) \\
& = & \Phi \circ \hat{\sigma}_\ell(p(z_1,\cdots,z_{\ell+1},z_1^*,\cdots,z_{\ell+1}^*)).
\eean
Since such elements are dense in $\cla_\ell$, we get
\[
\tau(a)=\Phi \circ \hat{\sigma}_\ell(a), \qquad  a \in \cla_\ell.
\]
 Thus by (\ref{lift}) $\tau$ admits the completely positive 
lifting $\hat{\sigma}_\ell$ and the result follows.
\qed

Thus one now has the following commutative diagram:

\hspace*{4em}
\def\labelstyle{\scriptstyle}
\xymatrix@C=25pt@R=40pt{
 & J_{\ell+1}\ar[r]\ar@{<->}[d]^{=} & A_{\ell+1}\ar[r]^{\sigma_\ell}\ar@{<->}[d]^{\psi_{\ell+1}} & 
         A_\ell \ar[dl]^{\hat{\sigma}_\ell} \ar[d]^{\tilde{\sigma}_\ell}&  \\
&E_{\ell+1}\ar[d]^{\seq}&  \psi_{\ell+1}(A_{\ell+1})\ar[d]^{\seq} & \cll(\clh_\ell\otimes C(S^1))\ar@{<->}[d]^{\cong}& \\
& M(E_{\ell+1}) \ar@{<->}[r]^{\cong} & \cll(F_{\ell+1}) &  \cll(F_{\ell+1}\oplus F_{\ell+1})\ar[l]_{C_\ell}
}

 \vspace{2ex}

Let $ev_1$ denote the following representation
of $C(S^1)$ on $\bbc$:
\[
ev_1(f)=f(1).
\]
Now take the trivial grading on $\bbc$. Then
$(\bbc,ev_1,0)$ gives an even Fredholm module for $C(S^1)$.

\blmma
The Fredholm module $(\bbc,ev_1,0)$ is a generator for the group $KK^0(C(S^1),\bbc)$.
\elmma
\prf
This can be seen as follows.
The identity projection gives a generating element
for $KK^0(\bbc,C(S^1))=K_0(C(S^1))=\bbz$.
The pairing of this with $[(\bbc,ev_1,0)]$ gives 1.
One can conclude from this that $[(\bbc,ev_1,0)]$ must be $\pm 1$.
\qed

\bppsn
$(\clh_\ell,\pi, \sgn D_{torus})]=
  (\clh_\ell\otimes C(S^1), \tilde{\sigma}_\ell, 2Q_\ell-I)\otimes_{ev_1} (\bbc,ev_1,0)$.
\eppsn
\prf
For this, one needs to note that
$(\clh_\ell\otimes C(S^1))\otimes\bbc\cong\clh_\ell$ where the tensor product
is the internal tensor product of Hilbert $C^*$-modules,
and
under this isomorphism, $(2Q_\ell-I)\otimes I$ is just the operator
$\sgn D_{torus}$.
\qed

Thus on multiplying the even Fredholm module $(\bbc,ev_1,0)$ from the left
by the $KK$-element we just computed, one gets the
odd fredholm module corresponding to the
spectral triple $(\clh_\ell,\pi_\ell,D_{torus})$ we have constructed
in the last section.

%%%%%%%%%%%%%%%%%%%%%%%%%%
%%%  BIBLIOGRAPHY
%%%%%%%%%%%%%%%%%%%%%%%%%%

%%%%%%%%%%%%%%%%%%%%%%%%%%%%%%%
\noindent{\sc Partha Sarathi Chakraborty}
(\texttt{parthac@imsc.res.in})\\
         {\footnotesize  Institute of Mathematical Sciences, 
CIT Campus, Chennai--600\,113, INDIA}\\[1ex]
{\sc Arupkumar Pal} (\texttt{arup@isid.ac.in})\\
         {\footnotesize Indian Statistical
Institute, 7, SJSS Marg, New Delhi--110\,016, INDIA}

%%%%%%%%%%%%%%%%%%%%%%%%%%%%%%%%%%%%%%%%%%%%%%%%%%%%%%%%%%%
%%%%%%%%%%%%%%%%%%%%%%%%%%%%%%%%%%%%%%%%%%%%%%%%%%%%%%%%%%%
%%%%%%%%%%%%%%%%%%%%%%%%%%%%%%%%%%%%%%%%%%%%%%%%%%%%%%%%%%%
%%%%%%%%%%%%%%%%%%%%%%%%%%%%%%%%%%%%%%%%%%%%%%%%%%%%%%%%%%%

\end{document}